\documentclass[a4paper]{amsart}
\usepackage{a4}

\usepackage{graphicx}
\usepackage[mathscr]{eucal}
\usepackage{amssymb}
\usepackage{epstopdf}
\usepackage{pstricks}
\usepackage{xy}
\xyoption{all}
\newcommand{\lk}{\ell\mathit{k}}
\newcommand{\bord}{\partial}
\newcommand{\IR}{\mathbb{R}}
\newcommand{\I}{\mathbf{1}}

\newcommand{\ZZ}{\mathbb{Z}}
\newcommand{\CC}{\mathbb{C}}
\newcommand{\SU}{\mathrm{SU}(2)}
\newcommand{\SL}{\mathrm{SL}_2(\CC)}
\newcommand{\PSL}{\mathrm{PSL}_2(\CC)}
\newcommand{\SO}{\mathrm{SO}(3)}

\newcommand{\sll}{\mathfrak{sl}_2(\CC)}
\newcommand{\Reg}[1]{\mathcal{R}\mathrm{eg}(#1)}

\newcommand{\im}{\mathop{\mathrm{im}}\nolimits}

\newcommand{\tr}{\mathop{\mathrm{tr}}\nolimits}
\newcommand{\myBP}{\mathscr{A}}
\newcommand{\myCS}{cs}
\newcommand{\myD}{D}
\newcommand{\myPr}{\pi}

\newcommand{\myT}{\partial M_K}
\newcommand{\myTr}{t}

\newcommand{\myS}{\omega}
\begin{document}
\theoremstyle{plain}
\newtheorem{theorem}{Theorem}
\newtheorem*{theorem*}{Theorem}
\newtheorem*{Main}{Main Theorem}
\newtheorem*{MultLemma}{Multiplicativity Lemma}
\newtheorem*{MainTheorem}{Main Theorem}
\newtheorem{prop}[theorem]{Proposition}
\newtheorem{corollary}[theorem]{Corollary}
\newtheorem{lemma}[theorem]{Lemma}
\newtheorem{claim}[theorem]{Claim}
\newtheorem*{conjecture*}{Conjecture}
\newtheorem{fact}[theorem]{Fact}
\theoremstyle{definition}
\newtheorem{definition}{Definition}
\newtheorem*{definition*}{Definition}
\newtheorem{example}{Example}
\newtheorem*{example*}{Example}
\newtheorem{notation}{Notation}
\newtheorem*{notation*}{Notation}
\newtheorem*{convention*}{Convention}
\theoremstyle{remark}
\newtheorem{remark}{Remark}
\newtheorem{question}{Question}

\title{On the asymptotic expansion of the colored Jones polynomial for torus knots} 
\author{J\'er\^ome Dubois \and Rinat Kashaev} 
\date{\today}                                           
\address{Centre de Recerca Matem\`atica \\ Apartat 50 \\ E--08193 Bellatterra  (Spain) \and Section de Ma\-th\'e\-ma\-ti\-ques \\ Universit\'e de Gen\`eve CP 64,
 2--4 Rue du Li\`evre \\ CH--1211 Gen\`eve 4 (Switzerland)}
\email{JDubois@crm.cat \and Rinat.Kashaev@math.unige.ch}
\begin{abstract}
In the asymptotic expansion of the hyperbolic specification of the colored Jones polynomial of torus knots, we identify different geometric
contributions, in particular Chern--Simons invaraint and Reidemeister torsion.  
\end{abstract}

\maketitle

\markboth{J\'er\^ome Dubois and Rinat Kashaev}{On the asymptotic expansion of the colored Jones polynomial for torus knots}
\pagestyle{myheadings}

\section{Introduction}

In this paper we consider the ``quantum hyperbolic invariant'' of a knot $K$ defined by the formula
\[
\langle K \rangle_N = \lim_{h \to 2\pi i/N} J'_N(K;h),\quad \forall
N\in\mathbb{Z}_{>1}.
\]
Here
\[
J'_N(K;h)=\frac{J_N(K;h)}{J_N(\bigcirc; h)}
\]
where $J_N(K; h)$ is the $N$-th colored Jones polynomial and
$\bigcirc$ stands for the unknot.  In the standard normalization 
one has
\[
J_N(\bigcirc; h) = \frac{\sinh(Nh/2)}{\sinh(h/2)}
\]
and it is known that in this normalization  the $N$-th colored Jones
polynomial vanishes at the point $h=2\pi i/N$ for any knot or
link. Thus, the quantity
$\langle K \rangle_N$ is a well--defined invariant.  

The aim of the present paper is to give a geometrical interpretation for the terms in the asymptotic expansion of  $\langle K \rangle_N$ at large $N$ in the case of torus knots.   This work is motivated by the ``volume conjecture" of \cite{Kas95,Mur} for the invariant $\langle K \rangle_N$ which states that
\begin{equation}
2 \pi \lim_{N \to \infty}\frac{\log | \langle K \rangle_N |}{
N} = v_3 |\!|S^3 \setminus K|\!|
\end{equation}
where $v_3$ denotes the 
hyperbolic volume of a regular ideal tetrahedron in
$\mathbb{H}^3$ and  $|\!| S^3 \setminus K |\!|$ is the simplicial or Gromov norm of the 3-manifold $S^3\setminus K$.

For any torus knot 
the asymptotic expansion of $\langle K \rangle_N$
 at large $N$ is derived in \cite{KasTir} to all orders, and one has $| \langle K \rangle_N | = O(N^{3/2})$. This is in fact the optimal estimation in the sense that there exists a sequence of integers of the form $N_j = 2pq(1+ 2 j)$ such that $\lim_{j \to \infty} | \langle K \rangle_{N_j} | / {N_j}^{3/2}$ exists and is not zero (see also~\cite{Zheng}, we thank H. Zheng for posing this question to us).
In particular, the volume conjecture appears to be trivially true in
 this case as  
 all torus knots are known to have vanishing simplicial norm.
 
In this paper, we prove that the non--abelian Reidemeister torsion and the Chern--Simons invariant appear in the asymptotic expansion of $\langle K \rangle_N$. The precise form of this result is as follows.

 For a knot $K \subset S^3$, let $N(K)$ be a tubular neighborhood of $K$, $M_K = S^3 \setminus N(K)$ its exterior. Let $\mathbb{T}^{K}_\lambda(\rho) \in \CC$ denote the sign--determined $\SL$-twisted Reidemeister
	torsion at an irreducible representation $\rho : \pi_1(M_K) \to \SL$ (with respect to the longitude $\lambda$ of the knot
	$K$) associated to $K$ defined by the first author
	in~\cite{JDFibre}, see Subsection~\ref{Torsion} of the present paper.  
Let $\Delta_{K}(t)$ denote the Alexander polynomial of a knot $K \subset S^3$ normalized so that $\Delta_K(t) = \Delta_K(t^{-1})$. According to 
the works of Milnor and Turaev~\cite{Milnor:1962,Turaev:1986}, the analytic map
\begin{equation}\label{TorsionAb}
\tau^{}_{K}(z) = \frac{2 \sinh(z)}{\Delta_{K}(e^{2z})}, \quad z \in \CC \setminus \{z \; |\; \Delta_{K}(e^{2z}) = 0 \},
\end{equation}
is essentially equal to the abelian Reidemeister torsion associated to $K$ (see~Section~\ref{AbTor} and especially Proposition~\ref{P:TorsionAlexander}). The poles  of $\tau^{}_K(z)$ describe the so--called \emph{bifurcation points} of the $\SL$-character variety of the knot exterior, corresponding to those abelian characters which are limits of non--abelian ones.

\begin{figure}[!tbh]
\begin{center}
\begin{pspicture}(4,7)
\psline[linewidth=.005](-.1,3)(4.1,3)
\psline[linewidth=.005](0,0)(0,6.5)\psline[linewidth=.005](4,0)(4,6.5)
\psdots[dotstyle=*](1.5,3)(3,3)
\psdots[dotstyle=*](0,3)(4,3)
\psdots[dotstyle=*](0,4.5)(0,6)
\psdots[dotstyle=*](4,2)(4,.5)
\psline(0,4.5)(4,.5)\psline(0,6)(4,2)
\psline(3,3)(4,3)\psline(4,3)(4,2)
\psline(0,6)(0,4.5)\psline(1.5,3)(3,3)
\rput[B]{0}(2.9,2.6){$P_\ell^+$}
\rput[B]{0}(1.4,2.6){$P_\ell^-$}
\rput[B]{0}(-.35,2.9){$O$}
\rput[B]{0}(4.35,2.9){$O'$}
\rput[B]{0}(-.35,4.4){$Q_\ell^-$}
\rput[B]{0}(-.35,5.9){$Q_\ell^+$}
\rput[B]{0}(4.35,.4){$R_\ell^-$}
\rput[B]{0}(4.35,1.9){$R_\ell^+$}
\rput[B]{0}(1.25,4){$A^{\diamond}_\ell$}
\rput[B]{0}(3.725,2.6){$A^{\triangleright}_\ell$}
\end{pspicture}
\end{center}
\caption{Euclidian picture associated to the $\ell$-th component of $X^{\mathrm{nab}}(M_K)$.}
\label{fig:1}
\end{figure}
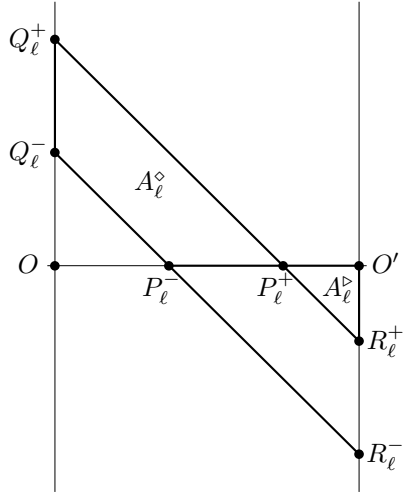

Let $K$ be the torus knot of type $(p,q)$. The non--abelian part of its
character variety, denoted by $X^{\mathrm{nab}}(M_K)$, consists of $N_{p,q} = (p-1)(q-1)/2$ connected
components, and each component intersects the abelian part in exactly
two bifurcation points. One can characterize the $\ell$-th non--abelian
component  by a unique pair of distinct positive integers $0<k^-_\ell<k^+_\ell<pq$
 satisfying certain conditions (see
Subsection~\ref{ModuliTorusK} for details and in particular Theorem~\ref{ChVTK}). In this case, the non--abelian Reidemeister torsion $\mathbb{T}^{K}_\lambda$ is locally constant on the $\SL$-character variety  (Proposition~\ref{propEx}). This is due to the structure of torus knot exteriors which are Seifert fibered spaces. We consider the euclidean
picture in Fig.~\ref{fig:1} where $O=(0,0)$, $O'=(\frac12,0)$,
$P_\ell^\pm=(\frac{k^\pm_\ell}{2pq},0)$, and two parallel segments passing through the points $P_\ell^\pm$ with slope $-pq$. In fact, this picture essentially
describes the real slice of $\ell$-th component of the character variety. 
Let $A^{\diamond}_\ell$ and
$A^{\triangleright}_\ell$ be the double areas of the
 trapezoid $P_\ell^-P_\ell^+Q_\ell^+Q_\ell^-$ and the triangle $P_\ell^+R_\ell^+O'$, respectively (a priori these are positive real numbers, but $A^{\diamond}_\ell$ is in fact an integer as we will see). We
 note that the set $\{A^{\triangleright}_\ell \pmod\ZZ\ \vert\
 1\le\ell\le N_{p,q}\}$ is essentially the
 Chern--Simons invariant of the knot $K$ (see Subsection~\ref{CSS} for
 details). Define the quantity

\begin{equation}\label{DefZ}
Z_N(K) = \sum_{\ell = 1}^{N_{p,q}} \varepsilon_\ell\sqrt{\left|\mathbb{T}^K_\lambda(\chi_\ell)\right|}\,
  A^{\diamond}_\ell  \, e^{- 2\pi i NA^{\triangleright}_\ell} ,\quad
\varepsilon_\ell=(-1)^{\left[ k^+_\ell/p\right]+
  \left[k^+_\ell/q\right]}
\end{equation}
where $N_{p,q} = {(p-1)(q-1)}/{2}$, $[x]$ is the integer part of the real number $x$, $\chi_\ell$ is any
 character in the $\ell$-th component (as the torsion is locally
 constant), and the positive value of the 
 square root is assumed. 

\begin{MainTheorem}\label{Main1}
For the $(p,q)$ torus knot $K$ the following asymptotic  equality holds: 
\begin{equation}
\label{MainEq}
2i^{pqN} {e^{\frac{i \pi}{2N}(\frac{p}{q} + \frac{q}{p})}} \cdot
\langle K \rangle_N = \sqrt{2pq} e^{i \pi/4}{N}^{3/2} Z_N(K) 
	 + (-1)^{pqN} \sum_{n=1}^\infty \frac{a_n(K)}{2^n n!} \left( \frac{i\pi}{pqN}\right)^{n-1}
\end{equation}
where $a_n(K)$ are finite type invariants.
\end{MainTheorem}

\begin{remark}
The sequence $a_n(K)$ is not only defined for torus knots but for all knots by the following general formula:
$$a_n(K) =\left. \frac{\partial^{2n} (z\tau_K(z))}{\partial z^{2n}} \right |_{z = 0}.$$
Each $a_n(K)$ is a finite type invariant (see Section~\ref{Proof1}).
\end{remark}
 
Equation~(\ref{MainEq}) gives a geometrical interpretation for each 
contribution in the asymptotic expansion of $\langle K \rangle_N$ 
in terms of classical invariants of knots. We can note
consistency with the \emph{Asymptotic Expansion Conjecture} 
(see~\cite[\S~7.2]{Ohtsuki}), namely appearances of the Chern--Simons 
invariant and the Reidemeister torsion, but we can also see some  
peculiarity in our formula, for example, the presence of the 
multiplication factors $A^{\diamond}_\ell$ and $\varepsilon_\ell$. 

We remark that, in~\cite{Mur2}, H. Murakami studied the asymptotic behavior of the colored Jones polynomial of torus knots 
$\log(J'_N\left(K;2\pi i r/N\right))/N$ in the complementary case $r \ne 1$. It would be interesting to understand the geometrical contributions in that case also.  Besides, in~\cite{Hikami} K. Hikami interprets the asymptotic expansion of the invariant $\langle K \rangle_N$ for $(2,q)$ torus knots from a  different point of view, namely in terms of  $q$-series identities.

\subsection*{Organization} In Section~\ref{modspace} we recall some
well known facts about the character varieties and fix the notation.
 In Section~\ref{Reidemeister}, we review the abelian and the non--abelian Reidemeister
 torsions for knots as presented in~\cite{JDFourier}, state the result
 of the computation for torus knots (Proposition~\ref{propEx}) and Yamaguchi's Theorem~\ref{Res} on the behavior of the non--abelian Reidemeister torsion near a bifurcation point. Section~\ref{CS} deals with the Chern--Simons invariant for knot exteriors. We compute it for torus knots using the technique of~\cite{KKK}.  In
 Section~\ref{Kashaev}, we review the quantum hyperbolic invariant for
 torus knots and, in particular, the integral representation
 of~\cite{KasTir}. Section~\ref{Proof1} contains the proof of the Main Theorem; the proof combines all results given in the preceding sections, especially the computation of the residues of the abelian Reidemeister torsion and of the Chern--Simons invariant for torus knots.

 \subsection*{Acknowledgements}
	The authors thank J.E. Andersen, S. Baseilhac, K. Hikami, V. Huynh Quang, T. Le, V. Turaev, Y. Yamaguchi and C. Weber for helpful discussions.
	
	This work is supported in part by the \emph{Swiss National Science Foundation}, the first author (J.D.) is also supported by the {European Community} with Marie Curie Intra--European Fellowship (MEIF--CT--2006--025316). While writing the paper, J.D. visited the CRM. He thanks the CRM for its hospitality.

\section{Review on $\SL$-character varieties}\label{modspace}

\subsection{Review on character variety}
\label{ModuliK}

Given a finitely generated group $\pi$ we let $$R(\pi; \SL) = \mathrm{Hom}(\pi; \SL)$$ denote the space of $\SL$-representations of $\pi$. This space is endowed with the {compact-open topology}. Here $\pi$ is assumed to have the discrete topology and the Lie group $\SL$ is endowed with the usual one. A representation $\rho \colon \pi \to \SL$ is called \emph{abelian} if $\rho(\pi)$ is an abelian subgroup of $\SL$. A representation $\rho$ is called \emph{reducible} if there exists a proper $U \subset \CC^2$ such that $\rho(g)(U) \subset U$, for all $g \in \pi$. Of course, any abelian representation is reducible. A non reducible representation is called \emph{irreducible}.

  The group $\SL$ acts on the representation space $R(\pi; \SL)$ by conjugation, but the naive quotient $R(\pi; \SL)/\SL$ is not Hausdorff in general. Following~\cite{CS:1983}, we will focus on the \emph{character variety} $X(\pi) = X(\pi; \SL)$ which is the set of \emph{characters} of $\pi$. Associated to $\rho \in R(\pi, \SL)$ is its character $\chi_\rho \colon \pi \to \CC$, defined by $\chi_\rho(g) = \mathrm{tr}(\rho(g))$. Let $g \in \pi$, following~\cite{CS:1983} we let $I_g \colon X(\pi) \to \CC$ denote the function defined by $I_g \colon \rho \mapsto \mathrm{tr}(\rho(g))$. Here $\mathrm{tr}(\rho(g))$ denotes the trace of the matrix $\rho(g)$.
  
  In some sense $X(\pi)$ is the ``algebraic quotient" of $R(\pi; \SL)$ by the action by conjugation of $\PSL$. It is well known that $R(\pi, \SL)$ and $X(\pi)$ have the structure of {complex algebraic affine sets} (see~\cite{CS:1983}).

  Let ${R}^\mathrm{irr}(\pi; \SL)$ denote the subset of irreducible representations of $\pi$ in $\SL$, and let $X^\mathrm{irr}(\pi)$ denote its image under the map ${R}(\pi; \SL) \to X(\pi)$. Note that two irreducible representations of $\pi$ in $\SL$ with the same character are conjugate by an element of $\SL$, see~\cite[Proposition 1.5.2]{CS:1983}.

For a knot $K$ in $S^3$, let $\Pi(K)$ denote its group. Let $\Pi(K)'$ be the subgroup generated by the commutators of $\Pi(K)$. It is well known that $\Pi(K)/\Pi(K)' \cong H_1(S^3 \setminus K; \ZZ) \cong \ZZ$ is generated by the meridian $\mu$ of $K$. As a consequence, each {abelian} representation of $\Pi({K})$ in $\SL$ is conjugate either to $\varphi_z \colon \Pi({K}) \to \SL$ defined by $\varphi_z(\mu) = \left(\begin{array}{cc}e^z & 0 \\0 & e^{-z}\end{array}\right)$, with $z \in \CC$, or to a representation $\rho$ with $\rho(\mu) = \pm \left(\begin{array}{cc}1 & 1 \\0 & 1\end{array}\right)$.  

Let $R^\mathrm{ab}(\Pi(K))$ denote the space of abelian representations, and let $X^\mathrm{ab}(\Pi(K))$ denote its image under the map $R(\Pi(K); \SL) \to X(\Pi(K))$. We write $X^\mathrm{nab}(\Pi(K))$ for the image of $R(\Pi(K); \SL) \setminus R^\mathrm{ab}(\Pi(K))$ under $R(\Pi(K); \SL) \to X(\Pi(K))$. Observe that $X^\mathrm{irr}(\Pi(K)) \subset X^\mathrm{nab}(\Pi(K))$ and that the converse is false.

If $M$ is a $3$-dimensional manifold, then we shall write $X^\mathrm{ab}(M) = X^\mathrm{ab}(\pi_1(M))$, $X^\mathrm{nab}(M) = X^\mathrm{nab}(\pi_1(M))$, $X^\mathrm{irr}(M) = X^\mathrm{irr}(\pi_1(M))$ etc.

\subsection{Character variety of torus knot groups}
\label{ModuliTorusK}

Let $\mathscr{T}(p,q)$ denote the torus knot of type $(p,q)$ where $p, q>1$ are coprime integers and set $\mathscr{M}(p,q) = M_{\mathscr{T}(p,q)}$.
 The group of $\mathscr{T}(p,q)$ admits the following well known presentation $$\Pi({p, q}) = \langle a, b \; |\; a^p = b^q \rangle.$$  Let $r, s \in \ZZ$ be such that $ps - qr = 1$. The meridian of $\mathscr{T}(p,q)$ is represented by the word $\mu = a^{-r} b^s$ and the longitude by $\lambda = a^p \mu^{-pq} = b^q \mu^{-pq}$. 

\begin{theorem}\label{ChVTK}
The non--abelian part $X^\mathrm{nab}(\mathscr{M}({p, q}))$ of the character variety $X(\mathscr{M}({p, q}))$ is the disjoint union of $N_{p,q} = (p-1)(q-1)/{2}$ copies of $\CC$ indexed by the pair $(\alpha, \beta)$ which satisfies the following conditions:
\begin{equation}\label{Conditions}
1 \leqslant \alpha \leqslant p-1, \ 1 \leqslant \beta \leqslant q-1, \ \alpha \equiv \beta \, (\mathrm{mod}\, 2).
\end{equation} 
On the component indexed by $(\alpha, \beta)$, we have $\chi(a) = 2 \cos \left( \frac{\alpha \pi}{p}\right)$, $\chi(b) = 2 \cos \left( \frac{\beta \pi}{q}\right)$. 

Moreover, we have $X^\mathrm{irr}(\mathscr{M}({p, q})) = X^\mathrm{nab}(\mathscr{M}({p, q})) \setminus X^{\mathrm{ab}}(\mathscr{M}({p, q}))$ and the component of $X^\mathrm{nab}(\mathscr{M}({p, q})) $ indexed by $(\alpha, \beta)$ is attached to $X^\mathrm{ab}(\mathscr{M}({p, q})) $ at two abelian representations $\varphi_{i\pi\frac{k^-}{pq}}$ and $\varphi_{i\pi\frac{k^+}{pq}}$ where $k^-$ and $k^+$ satisfy the following conditions:
\begin{equation}\label{EQkk}
0 < k^- < k^+ < pq,
\end{equation}
\begin{equation}\label{EQk}
k^- = \varepsilon^- (\beta ps - \varepsilon \alpha qr) \, (\mathrm{mod}\, pq),
\end{equation}
\begin{equation}\label{EQk'}
k^+ = \varepsilon^+ (\beta ps + \varepsilon \alpha qr) \, (\mathrm{mod}\, pq).
\end{equation}
Here $\varepsilon, \varepsilon^-, \varepsilon^+ \in \{\pm 1\}$.
\end{theorem}

We will use this description in the proof of  the Main Theorem. The first part of Theorem~\ref{ChVTK} is due to Klassen~\cite{Klassen:1991} and Le~\cite{Le}. Further notice the following remarks. 

\begin{remark}\label{remkk}
	The integers $k^-$ and $k^+$ satisfying conditions~(\ref{EQk}) and~(\ref{EQk'}) verify the following properties: 
	\begin{enumerate}
  \item $k^\pm$ is not divisible neither by $p$ nor by $q$;
  \item $k^+ \pm k^- \in 2\ZZ$;
  \item either $p$ divides $k^+ + k^-$, and $q$ divides $k^+-k^-$; or
 $q$ divides $k^+ + k^-$, and $p$ divides $k^+-k^-$.
\end{enumerate}
\end{remark}

\begin{remark}\label{remkk2} 
The quantities $m = \frac{k^+ - k^-}{2}$ and $l =
\frac{(k^+)^2 - (k^-)^2}{4pq}$ have the following properties:
\begin{enumerate}
  \item $m, l \in \ZZ_{>0}$,
  \item $m < \frac{pq}{2}$,
  \item either $p$ divides $m$ or $q$ divides $m$,
  \item $m$ divides $pq l$,
  \item $\frac{m^2}{pq} < l < m - \frac{m^2}{pq}$.
\end{enumerate}
\end{remark}

\begin{remark}
	In the case of $(2, q)$ torus knots  one has an explicit description of the character variety $X(\mathscr{M}(2, q))$, see~\cite{Le} for details. We use the fact that $\mathscr{T}(2, q)$ is a two-bridge knot and its group admits the following Wirtinger presentation:
\[
\Pi(2,q) = \langle a, b \; |\; wa = bw\rangle \text{ where } w = (ab)^{(q-1)/2}.
\]

	We write $x = I_a \colon \rho \mapsto \tr \rho(a)$, $y = I_{ab} \colon \rho \mapsto \tr \rho(ab)$. The character
	variety is parametrized by $x$ and $y$. Specifically, its
	abelian part $X^\mathrm{ab}(\mathscr{M}(p, q))$ is described
	by the equation $y - x^2 + 2 = 0$. The non--abelian part
	$X^\mathrm{nab}(\mathscr{M}(2, q))$ is described by the
	polynomial equation $\Phi(2, q) = 0$, where $\Phi(2, q)  \in
	\ZZ[x, y]$ satisfies the linear recurrence relation
\begin{equation}\label{ChV}
\Phi(2, q) = y \Phi(2, q-2) - \Phi(2, q-4) \text{ and } \Phi(2, 1) =
1, \Phi(2, 3) = y-1. 
\end{equation}
The non--abelian part $X^\mathrm{nab}(\mathscr{M}({2, q}))$ of the character variety is the disjoint union of $(q-1)/{2}$ copies of $\CC$ parametrized by $x$, indexed by $\ell \in \{1, \ldots, (q-1)/2\}$ and attached to $X^\mathrm{ab}(\mathscr{M}(2, q))$ at the $(q-1)$ abelian representations $\varphi_{z}$ for $$z \in \left\{i\pi  -i\pi\frac{2\ell - 1}{2q}, i\pi\frac{2\ell - 1}{2q}\right\}.$$
As a consequence, the $({q-1})/{2}$ pairs of integers $(k^-, k^+)$ are given by the following formulae: $$k^-_\ell = {2\ell - 1} \text{ and } k^+_\ell = 2q - 2\ell +1,$$
	where $\ell \in \left\{1, \ldots, ({q-1})/{2}\right\}$. Besides we have: 
\[
m_\ell = q - 2\ell + 1, \quad A_\ell^\diamond = \frac{q - 2\ell + 1}{2}, \quad 
A^\triangleright_\ell = \frac{(2\ell - 1)^2}{8q}.
\]
\end{remark}

\begin{remark}
	It is more complicated to compute the pairs $(k^-, k^+)$
	for arbitrary torus knots,  and in fact we do not know 
any ``closed'' formula for them. 
In Table~\ref{Table} we give some concrete examples.
	
\begin{table}[!h]
\begin{center}
\[\begin{array}{|c|c|c|c|c|c|}
\hline
(p,q)&N_{p,q}&(k^-_\ell,k^+_\ell)&m_\ell = \frac{k^+_\ell - k^-_\ell}{2}&\rule[-8pt]{0pt}{24pt} A_\ell^\diamond = \frac{{k^+_\ell}^2 - {k^-_\ell}^2}{4pq}&
A^\triangleright_\ell = \frac{{(pq - k^+_\ell)}^2}{4pq}\\[2pt]
\hline
\hline
(3,4)&3&(1,7)&3&1&\rule[-5pt]{0pt}{15pt}\frac{25}{48}\\
\cline{3-6}
&&(2,10)&4&2&\rule[-5pt]{0pt}{15pt}\frac{1}{12}\\
\cline{3-6}
&&(5,11)&3&2&\rule[-5pt]{0pt}{15pt}\frac{1}{48}\\
\hline
(3,5)&4&(1,11)&5&2&\rule[-5pt]{0pt}{15pt}\frac{4}{15}\\
\cline{3-6}
&&(2,8)&3&1&\rule[-5pt]{0pt}{15pt}\frac{49}{60}\\
\cline{3-6}
&&(4,14)&5&3&\rule[-5pt]{0pt}{15pt}\frac{1}{60}\\
\cline{3-6}
&&(7,13)&3&2&\rule[-5pt]{0pt}{15pt}\frac{1}{15}\\
\hline
(4,5)&6&(1,9)&4&1&\rule[-5pt]{0pt}{15pt}\frac{121}{80}\\
\cline{3-6}
&&(2,18)&8&4&\rule[-5pt]{0pt}{15pt}\frac{1}{20}\\
\cline{3-6}
&&(3,13)&5&2&\rule[-5pt]{0pt}{15pt}\frac{49}{80}\\
\cline{3-6}
&&(6,14)&4&2&\rule[-5pt]{0pt}{15pt}\frac{9}{20}\\
\cline{3-6}
&&(7,17)&5&3&\rule[-5pt]{0pt}{15pt}\frac{9}{80}\\
\cline{3-6}
&&(11,19)&4&3&\rule[-5pt]{0pt}{15pt}\frac{1}{80}\\
\hline
(4,7)&9&(1,15)&7&2&\rule[-5pt]{0pt}{15pt}\frac{169}{112}\\
\cline{3-6}
&&(2,26)&12&6&\rule[-5pt]{0pt}{15pt}\frac{1}{28}\\
\cline{3-6}
&&(3,11)&4&1&\rule[-5pt]{0pt}{15pt}\frac{289}{112}\\
\cline{3-6}
&&(5,19)&7&3&\rule[-5pt]{0pt}{15pt}\frac{81}{112}\\
\cline{3-6}
&&(6,22)&8&4&\rule[-5pt]{0pt}{15pt}\frac{9}{28}\\
\cline{3-6}
&&(9,23)&7&4&\rule[-5pt]{0pt}{15pt}\frac{25}{112}\\
\cline{3-6}
&&(10,18)&4&2&\rule[-5pt]{0pt}{15pt}\frac{25}{28}\\
\cline{3-6}
&&(13,27)&7&5&\rule[-5pt]{0pt}{15pt}\frac{1}{112}\\
\cline{3-6}
&&(17,25)&4&3&\rule[-5pt]{0pt}{15pt}\frac{9}{112}\\
\hline
\end{array}\]
\end{center}   
\caption{Numerical invariants of the non--abelian components of the character variety for some torus knots.}\label{Table}
\end{table}

\end{remark}

\section{Review on $\SL$-twisted Reidemeister torsion}
\label{Reidemeister}

\subsection{Preliminaries: sign-determined torsion of a CW--complex}

We review the basic notions and results about the sign-determined Reidemeister torsion introduced by Turaev which are needed in this paper. Details can be found in Milnor's survey~\cite{Milnor:1966} and in Turaev's monograph~\cite{Turaev:2002}.

\subsubsection*{Torsion of a chain complex}
Let $C_* = (\xymatrix@1@-.5pc{0 \ar[r] & C_n \ar[r]^-{d_n} & C_{n-1} \ar[r]^-{d_{n-1}} & \cdots \ar[r]^-{d_1} & C_0 \ar[r] & 0})$ be a chain complex of finite dimensional vector spaces over $\mathbb{C}$. Choose  a basis $\mathbf{c}^i$ for $C_i$ and  a basis $\mathbf{h}^i$ for the $i$-th homology group $H_i$. The torsion of $C_*$ with respect to these choice of bases is defined as follows.

Let $\mathbf{b}^i$ be a sequence of vectors in $C_{i}$ such that $d_{i}(\mathbf{b}^i)$ is a basis of $B_{i-1}= \im(d_{i} \colon C_{i} \to C_{i-1})$ and let $\widetilde{\mathbf{h}}^i$ denote a lift of $\mathbf{h}^i$ in $Z_i = \ker(d_{i} \colon C_i \to C_{i-1})$. The set of vectors $d_{i+1}(\mathbf{b}^{i+1})\widetilde{\mathbf{h}}^i\mathbf{b}^i$ is a basis of $C_i$. Let $[d_{i+1}(\mathbf{b}^{i+1})\widetilde{\mathbf{h}}^i\mathbf{b}^i/\mathbf{c}^i] \in \mathbb{R}^*$ denote the determinant of the transition matrix between those bases (the entries of this matrix are coordinates of vectors in $d_{i+1}(\mathbf{b}^{i+1})\widetilde{\mathbf{h}}^i\mathbf{b}^i$ with respect to $\mathbf{c}^i$). The \emph{sign-determined Reidemeister torsion} of $C_*$ (with respect to the bases $\mathbf{c}^*$ and $\mathbf{h}^*$) is the following alternating product (see~\cite[Definition 3.1]{Turaev:2000}):
\begin{equation}
\label{Def:RTorsion}
\mathrm{Tor}(C_*, \mathbf{c}^*, \mathbf{h}^*) = (-1)^{|C_*|} \cdot  \prod_{i=0}^n [d_{i+1}(\mathbf{b}^{i+1})\widetilde{\mathbf{h}}^i\mathbf{b}^i/\mathbf{c}^i]^{(-1)^{i+1}} \in \mathbb{C}^*.
\end{equation}
Here  $|C_*| = \sum_{k\geqslant 0} \alpha_k(C_*) \beta_k(C_*),$ where $\alpha_i(C_*) = \sum_{k=0}^i \dim C_k$ and  $\beta_i(C_*)  = \sum_{k=0}^i \dim H_k$.

The torsion $\mathrm{Tor}(C_*, \mathbf{c}^*, \mathbf{h}^*)$ does not depend on the choices of $\mathbf{b}^i$ and $\widetilde{\mathbf{h}}^i$. Further observe that if $C_*$ is acyclic (i.e. if $H_i = 0$ for all $i$), then $|C_*| = 0$.


\subsubsection*{Torsion of a CW--complex}
Let $W$ be a finite CW--complex and $\rho \in \mathrm{Hom}(\pi_1(W); \SL)$ a representation. We define the $\sll_{\rho}$-twisted cochain complex of $W$ to be
\[
C^*(W; \sll_\rho) = \mathrm{Hom}_{\pi_1(X)}(C_*(\widetilde{W}; \ZZ); \sll_\rho).
\]
Here $C_*(\widetilde{W}; \ZZ)$ is the complex of the universal covering with integer coefficients which is in fact a $\ZZ[\pi_1(W)]$-module (via the action of $\pi_1(W)$ on $\widetilde{W}$ as the covering group), and $\sll_\rho$ denotes the $\ZZ[\pi_1(W)]$-module via the composition $Ad \circ \rho$, where $Ad \colon \SL \to \mathrm{Aut}(\sll), A \mapsto Ad_A$ is the adjoint representation. This cochain complex $C^*(W; \sll_\rho)$ computes the {$\sll_\rho$-twisted cohomology} of $W$ which we denote as $H^*_\rho(W)$.

Let $\{e^{(i)}_1, \ldots, e^{(i)}_{n_i}\}$ be the set of $i$-dimensional cells of $W$. We lift them to the universal covering and we choose an arbitrary order and an arbitrary orientation for the cells $\left\{ {\tilde{e}^{(i)}_1, \ldots, \tilde{e}^{(i)}_{n_i}} \right\}$. If $\mathcal{B} = \{\mathbf{a}, \mathbf{b}, \mathbf{c}\}$ is an orthonormal basis of $\sll$, then we consider the corresponding ``dual" basis over $\CC$
$$\mathbf{c}^{i}_{\mathcal{B}} = \left\{ \tilde{e}^{(i)}_{1, \mathbf{a}}, \tilde{e}^{(i)}_{1, \mathbf{b}}, \tilde{e}^{(i)}_{1, \mathbf{c}}, \ldots, \tilde{e}^{(i)}_{n_i, \mathbf{a}}, \tilde{e}^{(i)}_{n_i, \mathbf{b}}, \tilde{e}^{(i)}_{n_i, \mathbf{c}}\right\}$$ of $C^i(W; \mathfrak{sl}(2)_\rho) = \mathrm{Hom}_{\pi_1(X)}(C_*(\widetilde{W}; \ZZ); \sll_\rho)$. Now choosing for each $i$ a basis $\mathbf{h}^{i}$ for the twisted cohomology $H^i_\rho(W)$, we can compute $$\mathrm{Tor}(C^*(W; \sll_\rho), \mathbf{c}^*_{\mathcal{B}}, \mathbf{h}^{*}).$$

The cells $\{ \tilde{e}^{(i)}_j \}^{}_{0 \leqslant i \leqslant \dim W, 1 \leqslant j \leqslant n_i}$ are in one-to-one correspondence with the cells of $W$, their order and orientation induce an order and an orientation for the cells $\{ e^{(i)}_j \}^{}_{0 \leqslant i \leqslant \dim W, 1 \leqslant j \leqslant n_i}$. Again, corresponding to these choices, we get a basis $c^i$ over $\IR$ for $C^i(W; \IR)$. 

Choose a \emph{cohomology orientation} of $W$, which is an orientation of the real vector space $H^*(W; \IR) = \bigoplus_{i\geqslant 0} H^i(W; \IR)$. Let $\mathfrak{o}$ denote this chosen orientation. Provide each vector space $H^i(W; \IR)$ with a reference basis $h^i$ such that the basis $\left\{ {h^0, \ldots, h^{\dim W}} \right\}$ of $H^*(W; \IR)$ is {positively oriented} with respect to $\mathfrak{o}$. Compute the sign-determined Reidemeister torsion $\mathrm{Tor}(C^*(W; \IR), c^*, h^{*}) \in \mathbb{R}^*$ of the resulting based and cohomology based chain complex and consider its sign $$\tau_0 = \mathrm{sgn}\left(\mathrm{Tor}(C^*(W; \IR), c^*, h^{*})\right) \in \{\pm 1\}.$$  

We define the sign-determined $Ad \circ \rho$-twisted Reidemeister torsion of $W$ to be
\begin{equation}\label{EQ:TorsionRaff}
\mathrm{TOR}(W; Ad \circ \rho, \mathbf{h}^{*}, \mathfrak{o}) = \tau_0 \cdot \mathrm{Tor}(C^*(W; \sll_\rho), \mathbf{c}^*_{\mathcal{B}}, \mathbf{h}^{*}) \in \mathbb{C}^*.
\end{equation}
This definition only depends on the combinatorial class of $W$, the conjugacy class of $\rho$, the choice of $\mathbf{h}^{*}$ and the cohomology orientation $\mathfrak{o}$. It is independent of the orthonormal basis $\mathcal{B}$ of $\sll$, of the choice of the lifts $\tilde{e}^{(i)}_j$, and of the choice of the positively oriented basis of $H^*(W; \IR)$. Moreover, it is independent of the order and the orientation of the cells (because they appear twice). 

One can prove that $\mathrm{TOR}$ is invariant under cellular subdivision, homeomorphism and simple homotopy equivalences. In fact, it is precisely the sign $(-1)^{|C_*|}$ in~(\ref{Def:RTorsion}) which ensures all these important invariance properties to hold.

\subsubsection*{Canonical orientation of knot exteriors}
In the case of knot exteriors in which we are interested, there exists a canonical cohomology orientation essentially defined by the meridian of the knot. The aim of this paragraph is to review it in details.

The exterior $M_K$ of $K$ is a $3$-dimensional CW--complex which has the same simple homotopy type as a $2$-dimensional CW--complex. We equip $M_K$ with its \emph{canonical cohomology orientation} defined as follows (see~\cite[Section V.3]{Turaev:2002}). We have 
$$H^*(M_K; \IR) = H^0(M_K; \IR) \oplus H^1(M_K; \IR)$$ 
and we base this $\IR$-vector space with $\{ \lbrack \! \lbrack pt \rbrack \! \rbrack, m^*\}$. Here $\lbrack \! \lbrack pt \rbrack \! \rbrack$ is the cohomology class of a point, and $\mu^* \colon \mu \mapsto 1$ is the dual of the meridian $\mu$ of $K$. This reference basis of $H^*(M_K; \IR)$ induces the so-called canonical cohomology orientation of $M_K$. In the sequel, we let $\mathfrak{o}$ denote the canonical cohomology orientation of $M_K$.

\subsection{Regularity for representations}

In this subsection we briefly review two notions of regularity (see~\cite{JDFibre} and~\cite{Porti:1997}). In the sequel $K \subset S^3$ denotes an oriented knot. We let $\Pi(K) = \pi_1(M_K)$ denote its group. The meridian $\mu$ of $K$ is supposed to be oriented according to the rule $\lk(K, \mu) = +1$, while the longitude $\lambda$ is oriented according to the condition $\mathrm{int}(\mu, \lambda) = +1$. Here $\mathrm{int}(\cdot, \cdot)$ denotes the intersection form on $\bord M_K$. 

We say that $\rho \in R^\mathrm{irr}(\Pi(K); \SL)$ is \emph{regular} if $\dim H^1_\rho(M_K) = 1$. This notion is invariant by conjugation and thus it is well--defined for irreducible characters.

\begin{example}
For the torus knot $\mathscr{T}(p, q)$, one can prove that each irreducible representation of $\Pi(p, q)$ in $\SL$ is regular. 

In the case of the figure--eight knot, one can also prove that each irreducible representation of its group in $\SL$ is regular.
\end{example}

Observe that for a regular representation $\rho$, we have $\dim H^1_\rho(M_K) = 1$, $\dim H^2_\rho(M_K) = 1$ and $H^j_\rho(M_K) = 0$ for all $j \ne 1, 2$. 

Let $\gamma$ be a simple closed unoriented curve in $\bord M_K$. Among irreducible representations we focus on the $\gamma$-regular ones. We say that $\rho \in R^\mathrm{irr}(\Pi(K); \SL)$ is \emph{$\gamma$-regular}, if (see~\cite[Definition 3.21]{Porti:1997}):
\begin{enumerate}
  \item the inclusion $\alpha \colon \gamma \hookrightarrow M_K$ induces an \emph{injective} map $$\alpha^* \colon H^1(M_K; \sll_\rho) \to H^1(\gamma; \sll_\rho),$$
  \item if $\tr(\rho(\pi_1(\bord M_K))) \subset \{\pm 2\}$, then $\rho(\gamma) \ne \pm \I$.
\end{enumerate} 
It is easy to see that this notion is invariant by conjugation and that $\gamma$-regularity implies regularity (the converse is false). Thus, for $\chi \in X^{\mathrm{irr}}(M_K)$ the notion of $\gamma$-regularity is well--defined. 

\begin{example}\label{ExReg}
For the torus knot $\mathscr{T}(p, q)$, one can prove that each irreducible representation of $\Pi(p, q)$ in $\SL$ is $\mu$-regular and also $\lambda$-regular.
\end{example}

Here is an alternative formulation, see~\cite[Proposition 3]{JDFibre}. Fix a generator $P^\rho$ of $H^0_\rho(\bord M_K)$. We recall that $$H^0_\rho(\bord M_K) = \sll^{Ad\circ \rho(\pi_1 \bord M_K)} = \left\{ v \in \sll \; |\; \forall g \in \pi_1 \bord M_K \; Ad_{\rho(g)}(v) = v \right\}.$$ The inclusion $\alpha \colon \gamma \hookrightarrow M_K$ and the cup product  induce the linear form $f^\rho_{\gamma} \colon H^1_\rho(M_K) \to \mathbb{C}$. We explicitly have $$f^\rho_{{\gamma}}(v) = B_{\sll} \left( {P^\rho, v(\gamma)}\right), \text{ for all } v \in  H^1_\rho(M_K).$$ 

\begin{prop}[Proposition 3 of \cite{JDFibre}]
A representation $\rho \in R^\mathrm{irr}(\Pi(K); \SL)$ is $\gamma$-regular if and only if the linear form $f^\rho_{\gamma} \colon H^1_\rho(M_K) \to \mathbb{C}$ is an isomorphism.
\end{prop}

\begin{remark}
A regular representation $\rho$ is $\gamma$-regular if and only if the linear form $f_\gamma^\rho$ is non degenerated (i.e. $f_\gamma^\rho \neq 0$).
\end{remark}

\subsection{Review on abelian Reidemeister torsion for knot exteriors}
\label{AbTor}

	The aim of this subsection is to compute the Reidemeister torsion of the exterior of $K$ twisted by the adjoint representation associated to an abelian representation of $G_K$ in terms of the Alexander polynomial of $K$.
	
	Let $\varphi_z : G_K \to \SL$ be the abelian representation such that $$\varphi_z(\mu) = \left(\begin{array}{cc}e^z & 0 \\0 & e^{-z}\end{array}\right)$$ and suppose that $\varphi_z$ is not boundary--central: $\varphi_z(\bord M_K) \not \subset \{\pm \I\}$  (i.e. $z \ne \pi k i$). When $e^{2z}$ is {not} a zero of the Alexander polynomial $\Delta_K$ of $K$ we say that $\varphi_z$ is {regular}. In this case, one can prove, following Klassen's arguments (\cite[Theorem 19]{Klassen:1991}), that $H^i_{\varphi_z}(M_K) \cong H^i(M_K; \ZZ) \otimes \CC$, for all $i$. 
	
	Let $h^{(0)} = P^\rho \in \sll$ be a fixed generator of $H^0_\rho(\bord M_K)$; then  $H^0_{\varphi_z}(M_K)$ is generated by $h^{(0)}$ and $H^1_{\varphi_z}(M_K)$ is generated by $h^{(1)} = h^{(0)}_1 + \cdots + h^{(0)}_{2n}$, where $h^{(0)}_k$ is the vector in $\sll^{2n}$ of which all entries are zero except the one of index $k$ which is equal to $h^{(0)}$.
	
	With this notation and choices, we have (see~\cite[Theorem 4]{Milnor:1962}, \cite[Subsection 1.2]{Turaev:1986} and~\cite[Proposition 4.4]{JDFourier})
	
\begin{prop}\label{P:TorsionAlexander}
	Let $\varphi_z$ be a regular abelian representation which is not boundary--central. The $(Ad \circ \varphi_z)$-twisted Reidemeister torsion of $M_K$ calculated in the basis $\{ h^{(0)}, h^{(1)}\}$ of $H^*_{\varphi_z}(M_K)$ and with respect to the canonical cohomology orientation of $M_K$ satisfies
\begin{equation}
\label{EQ:polynomeAlexander}
\mathrm{TOR}\left({M_K; Ad \circ \varphi_z, \{h^{(0)}, h^{(1)}\}}, \mathfrak{o} \right) =  - \tau_K(z) \tau_K(-z).
\end{equation}
Here $\tau_K(z)$ is the analytic map defined in Equation~(\ref{TorsionAb}).
\end{prop}
\begin{remark} If $K$ is the trivial knot, then $$\mathrm{TOR}\left({M_K; Ad \circ \varphi_z, \{h^{(0)}, h^{(1)}\}, \mathfrak{o}}\right) = 4\sinh^2(z)$$ is the twisted Reidemeister torsion of the solid torus $M_K$.
\end{remark}

\subsection{Review on non--abelian Reidemeister torsion for knot exteriors}
\label{Torsion}
This subsection gives a review of the constructions made in~\cite[\S~6]{JDFourier}. In particular, we shall explain how to construct distinguished bases for the twisted cohomology of knot exteriors.

\subsubsection*{How to construct natural bases for the twisted cohomology}

Let $\rho$ be a regular representation of $\Pi(K)$. One has a distinguished isomorphism induced by the cup product and the Killing form (which explicitly depends on the invariant vector $P^\rho$), see~\cite[Lemmas 5.1 \& 5.2]{JDFourier}: 
\[
\phi_{P^\rho} \colon H^2_\rho(M_K) \to H^2(M_K; \ZZ) \otimes \CC.
\]
Let $c$ be the generator of $H^2(\bord M_K; \ZZ) = \mathrm{Hom}(H_2(\bord M_K; \ZZ), \ZZ)$ 
corresponding to the fundamental class $\lbrack \! \lbrack \bord M_K \rbrack \! \rbrack \in H_2(\bord M_K; \ZZ)$ induced by the orientation of $\bord M_K$. The \emph{reference generator} of $H^2_\rho(M_K)$ is defined by 
\begin{equation}\label{EQ:Defh2}
h^{(2)}_\rho = \phi_{P^\rho}^{-1}(c).
\end{equation}

Let $\rho$ be a $\lambda$-regular representation of $\Pi(K)$. The \emph{reference generator} of $H^1_\rho(M_K)$ is defined by
\begin{equation}\label{EQ:Defh1}
h^{(1)}_\rho(\lambda) = (f^\rho_{\lambda})^{-1}(1).
\end{equation}

\subsubsection*{The Reidemeister torsion for knot exteriors}

Let $\rho \colon \Pi(K) \to \SL$ be a $\lambda$-regular representation. The \emph{Reidemeister torsion $\mathbb{T}^K_\lambda$} at $\rho$ is  
defined to be 
\[
\mathbb{T}^K_\lambda(\rho)  = \mathrm{TOR}\left( {M_K; Ad \circ \rho, \{h^{(1)}_\rho(\lambda), h^{(2)}_\rho\}, \mathfrak{o}} \right) \in \CC^*.
\]
It is an invariant of knots. Moreover, if $\rho_1$ and $\rho_2$ are two $\lambda$-regular representations which have the same character then $\mathbb{T}^K_\lambda(\rho_1) = \mathbb{T}^K_\lambda(\rho_2)$. Thus, $\mathbb{T}^K_\lambda$ defines a smooth map on the set $X^{\mathrm{irr}}_\lambda(M_K) = \{\chi \in X^{\mathrm{irr}}(M_K) \; |\; \chi \text{ is } \lambda\text{-regular}\} \subset \Reg{K}$.

\subsection{Reidemeister torsion for torus knots}
\label{RTTor}

The aim of this subsection is to state the computation of the $\SL$-twisted Reidemeister torsion for torus knots, see~\cite[\S~6.2]{JDFibre}. This result will be used to prove Theorem~\ref{Res}.

\begin{prop}[\cite{JDFibre}]\label{propEx}
If $\chi \in X^{\mathrm{irr}}(\mathscr{M}(p,q))$ lies in the component of the character variety indexed by the pair $(\alpha, \beta)$ which satisfies conditions~(\ref{Conditions}), then
\begin{equation}
\label{torsiontorus}
\mathbb{T}^{\mathscr{T}({p,q})}_{\lambda}(\chi) = \frac{16}{p^2q^2} \sin^2\left(\frac{\pi \alpha}{p}\right) \sin^2\left(\frac{\pi \beta}{q}\right).
\end{equation}
\end{prop}

\begin{remark}
	Example~\ref{ExReg} gives $X^{\mathrm{irr}}_\lambda(\mathscr{M}(p,q)) =X^{\mathrm{irr}}(\mathscr{M}(p,q)) $, thus $\mathbb{T}^{\mathscr{T}({p,q})}_{\lambda}$ is defined on the whole character variety of $\Pi(p,q)$ and not a priori only on $X^{\mathrm{irr}}_\lambda(\mathscr{M}(p,q)) $. Moreover, $\mathbb{T}^{\mathscr{T}({p,q})}_{\lambda}$ is locally constant on the character variety. 
\end{remark}

\begin{remark}
	The fact that $\mathbb{T}^{K}_{\lambda}$ is locally constant on the character variety in the case of torus knots is a particular property of such knots due to the fact that the exterior of a torus knot is a Seifert fibered manifold with two exceptional fibers (and a regular one). 

	In general, $\mathbb{T}^{K}_{\lambda}$ is, of course, not locally constant. An example is given by the figure--eight knot $4_1$, following~\cite{JDFibre}, we have
\begin{equation}\label{Torsionfigeight}
{(\mathbb{T}^{{4_1}}_\lambda(\rho))}^2 = \frac{1}{17 + 4 \, \mathrm{tr}(\rho(\lambda))}.
\end{equation}
\end{remark}

\subsection{Behavior of the non--abelian Reidemeister torsion near  a bifurcation point}

Let $z \in \CC$. Consider the abelian representation $\varphi_z\colon \Pi(K) \to \SL$ defined by:
\[
\varphi_z\colon \mu \mapsto \left(\begin{array}{cc}e^z & 0 \\0 & e^{-z}\end{array}\right)
\] 
and suppose that $\varphi_z$ is not boundary--central (i.e. $\varphi_z(\bord M_K) \not \subset \{\pm \I\}$). 

A result of G. Burde and G. de Rham~\cite{Burde,DeRham} states that there exists a reducible non--abelian representation $\rho_z \colon \pi_1(M_K) \to \SL$ which has the same character as $\varphi_z$ if and only if $\Delta_K(e^{2z}) = 0$. Furthermore, a recent result of Heusener, Porti and Su\'arez~\cite{HPS} states specifically that if $e^{2z}$ is a simple zero of $\Delta_K$ (i.e. $\Delta_K(e^{2z}) = 0$ and $\Delta'_K(e^{2z}) \ne 0$), then the corresponding representation $\rho_z$ is a \emph{bifurcation point}, i.e. $\rho_z$ is a limit of irreducible representations and is a smooth point of the $\SL$-representation variety (contained in a unique irreducible $4$-dimensional component of the $\SL$-representation variety). 

In the initial version of this paper, based on explicit calculations in the cases of torus knots and the figure--eight knot, we conjectured a relation between the residues of the analytic map $\tau_K$ and the limit values of the non--abelian Reidemeister torsion at the corresponding bifurcation points. Shortly after, this conjecture has been proved by Y. Yamaguchi in~\cite{YY2}. The result is as follows.

\begin{theorem}[\cite{YY2}, Theorem 1]\label{Res}
Let $z_0 \in \CC$ be such that $e^{2z_0}$ is a simple zero of the Alexander polynomial of $K$. Let $\chi_{z_0}$ denote the character corresponding to $\varphi_{z_0}$. The residue of $\tau^{}_{K}$ at $z_0$ satisfies:
\[
 \left(2 \underset{z = z_0}{\mathrm{Res}} \tau^{}_{K}(z)\right)^2_{} = \pm \lim_{\chi \to \chi_{z_0}} \mathbb{T}_\lambda^{K}(\chi).
\]
Here the limit is taken for irreducible $\SL$-characters $\chi$  which converge to $\chi_{z_0}$.
\end{theorem}

In the case of torus knots, which we are interested in that paper, Theorem~\ref{Res} is a trivial application of Formula~(\ref{torsiontorus}) and the proof reduces in this case to a direct computation of the residue at each bifurcation points.

\begin{proof}[Sketch of the proof of Theorem~\ref{Res}]
The proof is divided in two parts.
\begin{enumerate}
  \item{\emph{Existence part.}} The first point is to prove the existence of the limit:
\begin{equation}\label{lim}
\lim_{\underset{\chi \in X_\lambda^{\mathrm{irr}}(M_K)}{\mathrm{\chi \to \chi_{z_0}}}} \mathbb{T}_\lambda^{K}(\chi).
\end{equation}
The existence of limit~(\ref{lim}) is guaranteed by the fact that the reference bases of the twisted cohomology groups, defined in Equations~(\ref{EQ:Defh2}) and~(\ref{EQ:Defh1}), depend smoothly on the $\lambda$-regular character and can be smoothly extended  at $\chi_{z_0}$. Thus $\mathbb{T}_\lambda^{K}(\chi_{z_0})$ makes sense and is equal to limit~(\ref{lim}).
  \item{\emph{Computations part.}} With the preceding fact in mind, Theorem~\ref{Res} directly follows from the next lemma in which we compute the Reidemeister torsion at a reducible non--abelian representation.  
\begin{lemma}\label{LemmaMeta}
Let $\varphi$ be a reducible and non--abelian representation of $\Pi(K)$ in $\SL$ such that 
$$\varphi\colon \mu \mapsto \left(\begin{array}{cc}t & * \\0 & 1/t\end{array}\right), \; t \in \CC^*.$$
If we suppose that $t$ is a simple zero of $\Delta_K$ (i.e. $\Delta_K(t) = 0$ and $\Delta'_K(t) \ne 0$), then
\begin{equation}\label{TorsionMeta}
 \mathbb{T}_\lambda^{K}(\varphi) = \frac{(t-1)(1/t-1)}{\Delta'_K(t)\Delta'_K(1/t)}.
\end{equation}
\end{lemma} 
Essentially, Formula~(\ref{TorsionMeta}) can be considered as an analogue of Milnor--Turaev's formula for Reidemeister torsion at an abelian representation replacing abelian representations by appropriate reducible non--abelian ones and the Alexander polynomial by its derivative, see~\cite[Theorem 4]{Milnor:1962}, \cite[Subsection 1.2]{Turaev:1986}.
\end{enumerate}
\end{proof}

\section{$\SL$-Chern--Simons invariant for torus knots}\label{CS}

In the late 1980s, E. Witten considered a quantum field theory whose Lagrangian is the Chern--Simons functional. He argued that the Chern--Simons path integral on a $3$-manifold with an embedded link gives a (formal) $3$-dimensional interpretation of the Jones polynomial of links. Unfortunately path integral is not yet well--defined mathematically. But, the perturbative expansion method can be used and gives mathematical definition of knot invariants. Here we use another approach.

Let $M$ be a $3$-dimensional manifold. Let $\mathcal{A}$ denote the space of $\sll$-valued $1$-forms on $M$. The Chern--Simons functional $CS \colon \mathcal{A} \to \CC$ is defined by
\[
CS(A) = \frac{1}{8\pi^2} \int_M \mathrm{tr}(A \wedge dA + \frac{2}{3} A \wedge A \wedge A).
\]
In the case of a manifold with boundary this integral is not gauge invariant. Here we recall the definition of the Chern--Simons invariant for knot exteriors (and more generally for $3$-dimensional manifolds whose boundary consists of a $2$-dimensional torus) and explicitly compute it in the case of torus knots.

\subsection{Chern--Simons invariant for knot exteriors}
\label{SubsectionCS}

In this subsection we review the work of P.~Kirk and E.~Klassen~\cite{KKK} which allows us to compute the Chern--Simons invariant for torus knots.

Let $M$ denote a $3$-dimensional manifold whose boundary is non-empty and consists of a single $2$-dimensional torus $T = \bord M$. Following some ideas coming from Physics~\cite{RSW}, the Chern--Simons invariant of $M$ is considered in~\cite{KKK} as a lift of the restriction map induced by the inclusion $T \hookrightarrow M$:
\[
\xymatrix{ & \mathscr{B}_{T} \ar[d] \\ X(M) \ar[r] \ar[ur]^-{C_M} & X(T)}
\]
Here $\mathscr{B}_{T}$ is a $\CC^*$-bundle over $X(T)$ defined as follows. Let $(\mu, \lambda)$ be an oriented basis for $\pi_1(T)$. The map $t \colon X(T) \to \CC^3$ given by $\rho \mapsto (\tr \rho(\mu), \tr \rho(\lambda), \tr \rho(\mu \lambda))$ is an algebraic embedding of the character variety $X(T)$. Let $V(T)$ be the $2$-dimensional vector space (over $\CC$):
\[
V(T) = \mathrm{Hom}(\pi_1(T); \CC).
\]
The map $V(T) \to X(T)$ defined by $v \mapsto (x \mapsto 2\cos({2\pi iv(x)}))$ is a branched covering. The covering group $G$ is isomorphic to a semi--direct product of $\ZZ \oplus \ZZ$ and $\ZZ/2$ with the presentation
\[
G = \langle x, y, b \; |\; xyx^{-1}y^{-1} = bxbx = byby = b^2 = 1\rangle.
\]
Via the isomorphism $V(T) \to \CC^2$ defined by
\[
v \mapsto (v(\mu), v(\lambda))
\]
the action of $G$ on $V(T) \cong \CC^2$ is as follows
\[
x(\alpha, \beta) = (\alpha+1, \beta), \quad y(\alpha, \beta) = (\alpha, \beta +1), \quad b(\alpha, \beta) = (-\alpha, -\beta).
\]
Now we extend the action of $G$ to the the product $V(T) \times \CC^*$ by the formulas:
\begin{align}
\label{}
    x(\alpha, \beta; z) &= (\alpha + 1, \beta; ze^{2 \pi i \beta})   \\
    y(\alpha, \beta; z) &= (\alpha, \beta+ 1; ze^{-2 \pi i \alpha}) \\
    b(\alpha, \beta; z) &= (-\alpha, -\beta; z). 
\end{align}
Thus, the quotient $\CC^*$-bundle $\mathscr{B}_{T}$ over $X(T)$ is defined by the formula
\[
\mathscr{B}_T = V(T) \times \CC^*/G.
\]
We use the following notation for points in $\mathscr{B}_T$. Since $\mathscr{B}_T$ is a quotient of $\CC^2 \times \CC^*$ we write $[\alpha, \beta; z]$ for equivalence classes, so for example 
\[
[\alpha, \beta; z]=[\alpha + m, \beta + n; ze^{2\pi i(m\beta - n\alpha)}]. 
\]

The map $C_M \colon \rho \mapsto C_M(\rho) = [\gamma^{}_\mu, \gamma_\lambda^{}, e^{2\pi i CS_M(\rho)}]$ defines a lift of the restriction map induced by $T \hookrightarrow M$, see~\cite[Theorems 2.1 and 3.2]{KKK}. Here $(\gamma^{}_\mu, \gamma_\lambda^{})$ denotes a lift in $\CC^2$  of the restriction $\rho_{|\pi_1(T)}$. 

One can compute $C_M$ using the following Kirk--Klassen formula.

\begin{theorem}[\cite{KKK}, Theorem 2.7, Corollary 2.6 and Theorem 3.2]
\label{TKK}
Let $M$ denote an oriented $3$-dimensional manifold whose boundary $\bord M = T$ consist of a $2$-dimensional torus. Let $(\mu, \lambda)$ denote an oriented basis for $\pi_1(T)$.
\begin{enumerate}
  \item Let $\rho(t) \colon \pi_1(M) \to \SL$, $t \in [0, 1]$ be a path of representations. Let $(\gamma^{}_\mu(t), \gamma^{}_\lambda(t))$ denote a lift of $\rho(t)_{|\pi_1(T)}$ to $\CC^2$. Suppose
\[
C_M(\rho(t)) = [\gamma^{}_\mu(t), \gamma^{}_\lambda(t); z(t)]
\]
for all $t$.
Then 
\begin{equation}\label{C}
z(1) \cdot z(0)^{-1} = \mathrm{exp}\left( 2 \pi i \int_0^1 \gamma^{}_\mu(t) \gamma_\lambda'(t) - \gamma_\mu'(t) \gamma_\lambda^{}(t)\,dt\right).
\end{equation}
Furthermore, if $\rho(1)$ is the trivial representation, then $z(1) = 1$.
  \item There is an inner product 
  \begin{equation}\label{IP} 
  \langle \cdot, \cdot \rangle \colon \mathscr{B}_T \times \mathscr{B}_{-T} \to \CC^*
  \end{equation}
  given by taking the pair $\left( [\gamma^{}_\mu, \gamma_\lambda^{}, z], [\gamma^{}_\mu, \gamma_\lambda^{}, w]\right)$ to $z/w \in \CC^*$.
\end{enumerate}
\end{theorem}

\begin{remark}\label{orbi}
One can observe that Theorem~\ref{TKK} is true not only for $3$-dimensional manifolds but also for some $3$-dimensional orbifolds. For example, for the orbifold $W = D^2/\ZZ_2 \times S^1$ where the generator of the group $\ZZ_2$ acts on the $2$-dimensional disk $D^2$ by Euclidian rotation thru the angle $\pi$, the boundary $\bord W$ is the $2$-dimensional torus $T$. One has $\pi_1(W) = \ZZ_2 \times \ZZ$. If we identify the meridian $\mu$ and the longitude $\lambda$ of the boundary with the generators of $\ZZ_2$ and $\ZZ$ respectively,  then $X(W)$ is given by the equation 
\(
\chi(\mu) = \pm 2.
\) 
An easy application of Theorem~\ref{TKK} gives
\[
C_{W}(\chi) = \left[ n/2, \gamma_\lambda; e^{i \pi  n\gamma_\lambda}\right]
\]
for some integer $n$ and a character $\chi \in X(W)$, where for odd $n$ we have fixed the normalization by the condition $C_W(\chi_0) = [1/2, 0, 1]$ at the representation $\chi_0$ with trivial longitude. One can easily check that this is a well--defined lift of $X(W) \to X(T)$ to the bundle $\mathscr{B}_T$.
\end{remark}

\subsection{The areas $A^{\triangleright}_\ell$, $A^{\diamond}_\ell$ and the Chern--Simons invariant for knot exteriors}\label{CSS}

 The space $X(\myT)$ has the canonical
holomorphic simplectic form $\myS$, and
by using this additional structure we define a knot invariant as
follows.
 One has a branched 
covering mapping $\myPr\colon \CC^2\to X(\myT)$ given by the canonical
projection to the orbit space of the group action $(x,y)\mapsto (\pm
x+m,\pm y+n)$, $(m,n)\in\ZZ^2$.  One can choose a fundamental
domain for this action given by $\myD= D_1 \cap D_2$, where $D_1 = \{(x,y)\in\CC^2\vert\  0<\Re x<1/2\}$ and $D_2 = \{(x,y)\in\CC^2\vert\  |\Re y|<1/2\}$. Thus, topologically the set 
$X(\myT)$ is given by the quotient space $\overline{\myD}/\sim$ with
respect to the equivalence relation generated by $(i s,y)\sim (-i
s,-y)$, $(1/2+i s,y)\sim (1/2-is,-y)$, and $(x,1/2+i t)\sim
(x,-1/2+it)$, where $x,y\in\CC$, $s,t\in\IR$. The pullback of the holomorphic
symplectic form is given by the explicit formula
$\myPr^*\myS=2dx\wedge dy$.

Let $X_{\pm 2}(M_K)$ denote the finite sets of characters $\chi \in
X^\mathrm{irr}(M_K)$ such that $\chi(\mu) = \pm 2$, where $\mu$ is the
meridian of $K$.  Each of them is mapped to the set of branching points of $X(\bord M_K)$. One
can observe 
that $\sharp(X_{\pm 2}(\mathscr{T}(p,q))) = (p-1)(q-1)/2$ and if $K$ is
a hyperbolic knot then the characters corresponding to the lift (from the group $\mathrm{PSL}_2(\CC)$ to $\SL$) of the holonomy associated to the hyperbolic structure is in $X_{\pm 2}(M_K)$.  Thus in all
these cases, the sets $X_{\pm 2}(M_K)$ are non-empty.

Let $\myBP_K$ be the set of
  paths in the image of the non--abelian part of $X(M_K)$ in $X(\bord M_K)$  such that each $\gamma\in\myBP_K$
runs between the image of a bifurcation point  and that of a 
character in
$X_{\pm 2}(M_K)$ and admits a (unique) lift to a path $\tilde{\gamma}$ in $\overline{D_1}$. For any $\gamma\in\myBP_K$ let $\chi_\gamma\in X_{\pm 2}(M_K)$ be the end
point of $\gamma$. For each $\gamma\in\myBP_K$ we
associate a triangle $\myTr_\gamma$ in $\overline{D_1}$ as follows. If $\chi_\gamma \in X_2(M_K)$, then $\myTr_\gamma$ is the triangle with
sides $\tilde\gamma$ and two straight segments connecting the end
points of $\tilde\gamma$ 
with the point $(0,0)$. On the other hand, if $\chi_\gamma \in X_{-2}(M_K)$, then $\myTr_\gamma$ is the triangle with sides $\tilde\gamma$ and two straight segments connecting the end
points of $\tilde\gamma$ 
with the point $(1/2,0)$. See Fig.~\ref{fig:1} for the case of $\ell$-th component  of $X(\mathscr{M}(p, q)$). We define a function
$\myCS_K\colon \myBP_K\to\CC$ which associates to $\gamma\in \myBP_K$ the
complex symplectic area of the projected triangle
  $\myPr(\myTr_\gamma)$. This function generalizes the Chern--Simons invariant
  defined in~\cite{KKK}. Suppose that $K$ is the $(p, q)$ torus knot. Let $\gamma$ be such that $\chi_\gamma \in X_{-2}(M_K)$ and $\chi_\gamma$ is on the $\ell$-th component of $X(M_K)$, then $|cs_K(\gamma)| = A^{\triangleright}_\ell$.

	The exterior of the knot $M_K$ and the orbifold $W(K) = N(K) / \ZZ_2$, where $\ZZ_2$ acts by rotation on the first component of $N(K) = D^2 \times S^1$, have the same boundary. In $W(K)$ the meridian $\mu$ bounds a $2$-dimensional orbifold $D^2/\ZZ_2$, hence each character $\chi' \in X(W(K))$ satisfies $\chi'(\mu) = \pm 2$. Let $\chi \in X(M_K)$ be such that  $\chi$ and $\chi'$ coincide on the peripheral system of $K$, then $\chi \in X_{\pm 2}(M_K)$. 
 	We define on $X_{\pm 2}(M_K)$ a complex valued discrete function by
\begin{equation}\label{SCS}
CS_K \colon X_{\pm 2}(M_K) \to \CC^*, \quad CS_K(\chi) = \left\langle {C_{M_K}(\chi), C_{W(K)}(\chi') }\right\rangle,
\end{equation}
where $\chi$ and $\chi'$ coincide on the peripheral system of $K$. 

By construction $CS_K$ is a knot invariant. In this paper we call this discrete function the \emph{Chern--Simons invariant of the knot $K$}. Let $\gamma \in \mathscr{A}_K$ and let $\chi_\gamma$ be the end point of $\gamma$. We have $$CS_K(\chi_\gamma) = e^{2\pi i cs_K(\gamma)}.$$

Let $\gamma_1$, $\gamma_2 \in \myBP_K$ have coinciding end points, but start from two bifurcation points on one and the same connected components of $X(M_K)$. Then the two corresponding triangles have areas differing by an integer (since both of them give the same $CS_K(\chi_\gamma)$). In the case of torus knots the difference of those areas is precisely the symplectic area $A^{\diamond}_\ell$  of the trapezoid in Fig.~\ref{fig:1}.

\subsection{Computation of the Chern--Simons invariant for torus knots}
\label{CSTor}

\begin{prop}\label{Prop:CSKK}
Let $\rho : \Pi(p,q) \to \SL$ be a non--abelian representation whose character lies in the component of $X(\mathscr{M}(p,q))$ parametrized by $(\alpha, \beta)$ (see Theorem~\ref{ChVTK}). If the matrix $\rho(\mu)$ is conjugate to
\[
\rho(\mu) \simeq \left(\begin{array}{cc}e^{2\pi i \gamma_\mu} & * \\0 & e^{-2\pi i \gamma_\mu}\end{array}\right)
\]
then
\begin{equation}\label{CSKK}
C_{\mathscr{M}(p,q)}(\rho) = \left[ {\gamma_\mu, \frac{1}{2} - pq\gamma_\mu; \exp \left({2\pi i\left( \frac{(\beta ps + \varepsilon \alpha qr)^2}{4pq} - \frac{\gamma_\mu}{2}\right)}\right)} \right]. 
\end{equation}
Here $ps - qr =1$ and the result does not depend on the choice of $\varepsilon \in \{\pm 1\}$.
\end{prop}

The proof is based on Kirk--Klassen's Theorem~\ref{TKK}.

\begin{proof}[Proof of Proposition~\ref{Prop:CSKK}]
	Our computation is based on the fact that a non--abelian representation $\rho$ of $\Pi(p,q)$ which lies in the component parametrized by $(\alpha, \beta)$ is connected to the trivial representation $\vartheta$ by a path of representations which contains the bifurcation point $\varphi_{z_k}$, where $z_k = i\pi\frac{k}{pq}$ with $k$ satisfying conditions of Theorem~(\ref{ChVTK}). This path is divided into two distinct parts connected by $\varphi_{z_k}$:
\begin{enumerate}
  \item the path $(\varphi_t)_{0\leqslant t \leqslant \theta}$ of abelian representations which connects $\vartheta$ to $\varphi_{z_k}$,
  \item the path $(\rho_t)_{0\leqslant t \leqslant 1}$ of non--abelian representations which connects $\rho = \rho_0$ to $\varphi_{z_k} = \rho_1$.
\end{enumerate}

	The computation is done in two steps.
\begin{enumerate}
  \item First, applying the second part of Kirk--Klassen's Theorem~\ref{TKK} to the abelian representation $\varphi_{z_k}$ where $z_k = i\pi\frac{k}{pq}$, we get:
\begin{equation}\label{CSM}
C_{\mathscr{M}(p,q)}(\varphi_{z_k}) = \left[ {\frac{k}{2pq}, 0; 1} \right].
\end{equation}
Observe that for any abelian representation $\varphi$ one has $\varphi(\lambda) = \I$.
  \item Now, we apply Kirk--Klassen's Theorem~\ref{TKK} to the path $(\rho_t)_{0\leqslant t \leqslant 1}$ of non--abelian representations. Suppose that $C_{\mathscr{M}(p,q)}(\rho_t) = [\gamma_\mu(t), \gamma_\lambda(t); z(t)]$.
  As $\lambda = a^p \mu^{-pq}$, we have $$\gamma_\lambda(t) = \frac{1}{2} - pq\gamma_\mu(t) \; (\text{because } a^p \text{ is central in } \Pi(p,q) \text{ and thus } \rho_t(a^p) = -\textbf{1})$$ and thus $$\gamma_\lambda'(t) = -pq \gamma_\mu'(t).$$ Moreover we have $\gamma_\lambda(1) = \frac{1-k}{2}$. Choosing appropriate representative for Formula~(\ref{CSM}), one has 
  \[
  C_{\mathscr{M}(p,q)}(\varphi_{z_k}) = \left[ {\frac{k}{2pq}, 0; 1} \right] = \left[{\frac{k}{2pq}, \frac{1-k}{2}; \mathrm{exp}\left( -2\pi i \frac{k(1-k)}{2pq}\right)} \right].
  \]

With this, Formula~(\ref{C}) gives us
\[
z(1)\cdot z(0)^{-1} = \mathrm{exp}\left( -i\pi\int_0^1 \gamma'_\mu(t) \,dt \right) = \mathrm{exp}\left( 2\pi i \left( \frac{\gamma_\mu}{2} - \frac{k}{4pq}\right)\right).
\]  
Thus,
\[
z(0) = \mathrm{exp}\left( 2\pi i\left(  \frac{k^2}{4pq}  - \frac{\gamma_\mu}{2}\right)\right)
\]
which achieves the proof.
\end{enumerate}
\end{proof}

We are now ready to compute the discrete function $CS_{\mathscr{T}(p,q)} \colon X_{\pm 2}(\mathscr{M}(p, q)) \to \CC^*$ defined in equation~(\ref{SCS}). For each $1 \leqslant \ell \leqslant (p-1)(q-1)/2$, let $\chi_\ell^\pm$ be the unique irreducible character of $\Pi(p,q)$ with $\chi(\mu) = \pm 2$ and which lies in $\ell$-th connected component of $X^{\mathrm{irr}}(\mathscr{M}(p,q)$. We have $X_{\pm 2}(\mathscr{M}(p, q)) = \{\chi_\ell^\pm \ | \ 1 \leqslant \ell \leqslant (p-1)(q-1)/2\}$.

\begin{corollary}\label{Cor:CSKK}
The map $CS_{\mathscr{T}(p,q)} \colon X_{\pm 2}(\mathscr{M}(2, q)) \to \CC^*$ defined in equation~(\ref{SCS}) satisfies
\begin{equation}\label{CSKK2}
CS_{\mathscr{T}p,q)}(\chi^+_\ell) = e^{i\pi\frac{{k^-_\ell}^2}{2pq}} \text{ and } CS_{\mathscr{T}(p,q)}(\chi^-_\ell) = e^{- i\pi\frac{{(pq - k^+_\ell)}^2}{2pq}}.
\end{equation}
Here $k^\pm_\ell$ are the two integers which satisfies conditions of Subsection~\ref{ModuliTorusK}.
\end{corollary}

\begin{remark}
We have $$CS_{\mathscr{T}p,q)}(\chi^+_\ell) = e^{i\pi\frac{{k^+_\ell}^2}{2pq}}.$$
This is just the consequence of the fact that $A^{\diamond}_\ell = \frac{{k^+_\ell}^2 - {k^-_\ell}^2}{4pq}$ is an integer  (see Remark~\ref{remkk2}).
\end{remark}

\begin{remark}
The Chern--Simons invariant for torus knots has also been discussed  in
papers \cite{Hikami,HK}.
\end{remark}

\begin{proof}[Proof of Corollary~\ref{Cor:CSKK}]
The computation of $C_{W({\mathscr{T}(p, q)})}(\rho)$ for any character of $\Pi(p, q)$ in $\SL$ is done in Remark~\ref{orbi}. Now, applying the inner product~(\ref{IP}) we get
\[
CS_{\mathscr{T}(p,q)}(\chi^+_\ell) = \left\langle {C_{\mathscr{M}(p, q)}(\chi^+_\ell), C_{W({\mathscr{T}(p, q)})}({\chi^+_\ell}')}\right\rangle = e^{i\pi\frac{{k^-_\ell}^2}{2pq}}
\] 
and
\[
CS_{\mathscr{T}(p,q)}(\chi^-_\ell) = \left\langle {C_{\mathscr{M}(p, q)}(\chi^-_\ell), C_{W({\mathscr{T}(p, q)})}({\chi^-_\ell}')}\right\rangle = e^{- i\pi\frac{{(pq - k^+_\ell)}^2}{2pq}}.
\] 
Here $\chi_\ell^\pm$ and ${\chi_\ell^\pm}'$ coincide on the peripheral system of the knot.
\end{proof}

\section{Review on the quantum hyperbolic invariant for torus knots}\label{Kashaev}

Using the relationship between the colored Jones polynomial and the quantum hyperbolic invariant, one has the following result.
\begin{prop}[\cite{KasTir}, Lemma 2]\label{IntegralKashaev}
Let $p$ and $q$ be two coprime integers. For the torus knot $\mathscr{T}({p, q})$ of type $(p,q)$, the quantum hyperbolic invariant has the following  complex integral representation:
\begin{equation}\label{intkastir}
2 \cdot \langle \mathscr{T}({p, q}) \rangle_N = \left( \frac{pq N}{2}\right)^{3/2} \cdot e^{-\frac{i\pi}{2N}\left( \frac{p}{q} + \frac{q}{p} + \frac{N}{2}\right)} \cdot \int_{\mathscr{C}} e^{\pi pq N(z + \frac{i}{2}z^2)} z^2 \tau_{\mathscr{T}({p, q})}(z) \, dz.
\end{equation}
Here the path of integration $\mathscr{C}$ is the image of the real line under the mapping $\IR \ni x \mapsto xe^{i \phi} \in \CC$ where $\phi$ is to be chosen by the convergence condition.
\end{prop}

Proposition~\ref{IntegralKashaev} is the main ingredient to prove:

\begin{prop}[\cite{KasTir}, Theorem]\label{KasTirP}
Let $p$ and $q$ be two coprime integers. For the torus knot $\mathscr{T}({p, q})$ of type $(p,q)$, the quantum hyperbolic invariant has the following asymptotic expansion at large $N$:
\begin{equation}\label{KasTirEq}
e^{\frac{i\pi}{2N}(\frac{p}{q} + \frac{q}{q})} \cdot \langle \mathscr{T}(p,q) \rangle_N = \sum_{k = 1}^{pq-1} \langle \mathscr{T}(p,q) \rangle^{(k)}_N + \langle \mathscr{T}(p,q) \rangle^{(\infty)}_N.
\end{equation}
Here
\begin{equation}\label{Tkk}
\langle\mathscr{T}(p,q) \rangle^{(k)}_N =  \left( \frac{pq N}{2}\right)^{3/2} \cdot \frac{e^{-\frac{i\pi}{4}}}{2} \cdot \underset{z = i\pi\frac{k}{pq}}{\mathrm{Res}} \left({e^{\pi pq N(z + \frac{i}{2}z^2)} z^2 \tau_{\mathscr{T}({p, q})}(z)}\right),
\end{equation}
and
\begin{equation}\label{Tinf}
\langle\mathscr{T}(p,q) \rangle^{(\infty)}_N = \frac{i^{  pq N}}{4} \sum_{n=0}^\infty \frac{a_n(K)}{n!} \left( \frac{i\pi}{2pqN}\right)^{n-1}
\end{equation}
where $a_n(K) =\left. \frac{\partial^{2n} (z\tau_K(z))}{\partial z^{2n}} \right |_{z = 0}$.
\end{prop}

Formula~(\ref{KasTirEq}) is our starting point in the proof of the Main Theorem. Moreover a direct computation gives:
\begin{equation}\label{Tk}
\langle\mathscr{T}(p,q) \rangle^{(k)}_N = 2 \left(\frac{N}{2pq}\right)^{3/2} \cdot e^{\frac{i\pi}{4}} \cdot (-1)^{(N - 1)k} \cdot e^{-\frac{i\pi k^2}{2pq} N}\cdot k^2 \cdot \sin({\pi k}/{p})\sin({\pi k}/{q}).
\end{equation}

\section{Proof of the Main Theorem}
\label{Proof1}

The proof is a direct computation of each part of equality~(\ref{MainEq}) and combines Propositions~\ref{KasTirP} and~\ref{Prop:CSKK}.

\begin{proof}[Proof of the Main Theorem]

For each $\ell \in \{1, \ldots, (p-1)(q-1)/2\}$, the $\ell$-th connected component of $X^{\mathrm{irr}}(\mathscr{M}(p,q)$ intersects the abelian one at two bifurcation points $\varphi^{}_{i\pi\frac{k^+_\ell}{pq}}$ and $\varphi_{i\pi\frac{k^-_\ell}{pq}}$. Here $k^\pm_\ell $ are the two integers which satisfy conditions of Subsection~\ref{ModuliTorusK}. We rewrite the sum $\sum_{k = 1}^{pq-1} \langle \mathscr{T}(p,q) \rangle^{(k)}_N$, where $\langle \mathscr{T}(p,q) \rangle^{(k)}_N$ is defined in equation~(\ref{Tk}), as follows. 

First, observe that $\langle \mathscr{T}(p,q) \rangle^{(k)}_N = 0$ for all integers $k$ which are divisible by $p$ or $q$. Write
\[
\langle \mathscr{T}(p,q) \rangle^{(j)}_N = \Gamma_{j} \cdot \sin\left({\pi j}/{p}\right)\sin\left({\pi j}/{q}\right),
\]
where
\[
\Gamma_{j} = 2 \left(\frac{N}{2pq} \right)^{3/2}e^{\frac{i\pi}{4}}e^{-i\pi N \frac{j^2}{2pq}} (-1)^{(N - 1)j}j^2.
\]

Next observe that the integers $k^\pm_\ell$ indexed by $\ell \in \{1, \ldots, (p-1)(q-1)/2\}$ exhaust all the integers between $1$ and $pq - 1$ which are mutually prime with $pq$.
Thus, we  have
\[
\sum_{j = 1}^{pq-1} \langle \mathscr{T}(p,q) \rangle^{(j)}_N = \sum_{\ell = 1}^{{(p-1)(q-1)}/{2}} \langle \mathscr{T}(p,q) \rangle^{(k^-_\ell)}_N + \langle \mathscr{T}(p,q) \rangle^{(k^+_\ell)}_N.
\]
It is easy to observe that (see Remaks~\ref{remkk} \&~\ref{remkk2})
\[
\sin\left( {\pi k^+_\ell}/{p}\right)\sin\left( {\pi k^+_\ell}/{q}\right) = - \sin\left({\pi k^-_\ell}/{p}\right)\sin\left({\pi k^-_\ell}/{q}\right).
\]
As a consequence, we obtain
\[
\sum_{j = 1}^{pq-1} \langle \mathscr{T}(p,q) \rangle^{(j)}_N = \sum_{\ell = 1}^{(p-1)(q-1)/2} (\Gamma_{k^+_\ell}- \Gamma_{k^-_\ell}) \cdot \sin\left({\pi k^+_\ell}/{p}\right)\sin\left({\pi k^+_\ell}/{q}\right).
\]
Using the properties in Remark~\ref{remkk2} (1), we have
\[
\Gamma_{k^+_\ell}- \Gamma_{k^-_\ell} = 2 \left(\frac{N}{2pq} \right)^{3/2}e^{\frac{i\pi}{4}}e^{-i\pi N \frac{{k^+_\ell}^2}{2pq}} (-1)^{(N - 1)k^+_\ell}({k^+_\ell}^2 - {k^-_\ell}^2).
\]
Thus,
\[
\sum_{j = 1}^{pq-1} \langle \mathscr{T}(p,q) \rangle^{(j)}_N = {N^{3/2}} {e^{\frac{i\pi}{4}}}{\sqrt\frac{pq}{2}}   \cdot \Sigma
\]
where
\begin{equation}\label{EQS}
\Sigma = \sum_{\ell = 1}^{(p-1)(q-1)/2} \frac{(-1)^{(N- 1)k^+_\ell}}{{4pq}} {({k^+_\ell}^2 - {k^-_\ell}^2)} \cdot e^{-i\pi N \frac{{k^+_\ell}^2}{2pq}} \cdot \frac{4}{pq} \sin\left( {\pi k^+_\ell}/{p}\right)\sin\left({\pi k^+_\ell}/{q}\right)
\end{equation}

Here is a geometric explanation of the quantities which appear in equation~(\ref{EQS}). 

Using Theorem~\ref{Res}, one has
\[
 \sqrt{\left|\mathbb{T}^{\mathscr{T}(p,q)}_\lambda(\chi_\ell)\right|}= 
 \frac{4}{pq}\sin\left( {\pi k^+_\ell}/{p}\right)\sin\left({\pi k^+_\ell}/{q}\right) = (-1)^{k^+_\ell-1} 2\underset{z = i\pi\frac{k^+_\ell}{pq}}{\mathrm{Res}} \tau^{}_{\mathscr{T}(p,q)}(z).
\]

Corollary~\ref{Cor:CSKK} gives
\[
e^{i\pi\frac{{k^+
_\ell}^2}{2pq}} = {CS}_{\mathscr{T}(p,q)}({\chi^+_{\ell}})
\]
where $\chi^+_{\ell}$ lies in the $\ell$-th component  of $X^\mathrm{irr}(\mathscr{M}(p,q))$ and satisfies $\chi^+_\ell(\mu) = 2$ (see Corollary~\ref{Cor:CSKK}). Writing ${{CS}_\ell}$ instead of ${CS}_{\mathscr{T}(p,q)}({\chi^+_{\ell}})$ we arrive at the formula
\[
\sum_{j = 1}^{pq-1} \langle \mathscr{T}(p,q) \rangle^{(j)}_N = N^{3/2}{e^{\frac{i\pi}{4}}}{\sqrt{\frac{pq}{2}}}\; \cdot \; \Sigma'
\]
where $$\Sigma' = \sum_{\ell = 1}^{(p-1)(q-1)/2} \frac{(-1)^{(N - 1)k^+_\ell}}{4pq}({k^+_\ell}^2 - {k^-_\ell}^2) \cdot {CS}_\ell^{-N} \cdot  \sqrt{\left|\mathbb{T}^{\mathscr{T}(p,q)}_\lambda(\chi_\ell)\right|}. $$

Now we compute the symplectic areas $A^{\diamond}_\ell$ and
$A^{\triangleright}_\ell$. Recall that $A^{\diamond}_\ell$ and
$A^{\triangleright}_\ell$ are respectively the double areas of the
 trapezoid $P_\ell^-P_\ell^+Q_\ell^+Q_\ell^-$ and the triangle $P_\ell^+R_\ell^+O'$, see Fig.~\ref{fig:1}. So, one has
 
 \[
 A^{\diamond}_\ell = \frac{{k^+_\ell}^2 - {k^-_\ell}^2}{4pq} \text{ and } A^{\triangleright}_\ell = \frac{(pq - k_\ell^+)^2}{4pq}.
 \]
 
 Combining the two last formulas, we finally obtain
 \[
 \sum_{j = 1}^{pq-1} \langle \mathscr{T}(p,q) \rangle^{(j)}_N = {N}^{3/2} \sqrt{\frac{pq}{2}} \cdot \frac{e^{\frac{i\pi}{4}}}{i^{pq N}}  \sum_{\ell = 1}^{(p-1)(q-1)/2} \varepsilon_\ell \sqrt{\left|\mathbb{T}^{\mathscr{T}(p,q)}_\lambda(\chi_\ell)\right|} A^{\diamond}_\ell e^{- 2\pi i N A^{\triangleright}_\ell},
 \]
where $\varepsilon_\ell = (-1)^{[k_\ell^+/p] + [k_\ell^+/q]}$, thus
\[
 \sum_{j = 1}^{pq-1} \langle \mathscr{T}(p,q) \rangle^{(j)}_N = \sqrt{\frac{pq}{2}} \cdot \frac{e^{\frac{i\pi}{4}}}{i^{pq N}}\cdot {N}^{3/2} Z_N(\mathscr{T}(p,q))
\]
which achieves the proof of Formula~(\ref{MainEq}).

We now prove that $a_n(K)$ is a finite type invariant for each $n$.
One has
$$a_n(K) =\left. \frac{\partial^{2n} (z\tau_K(z))}{\partial z^{2n}} \right |_{z = 0}, \text{ where } \tau^{}_{K}(z) = \frac{2 \sinh(z)}{\Delta_{K}(e^{2z})}$$
As a consequence, $a_n(K)$ is a linear combination of products of coefficients of  the Alexander polynomial $\Delta_K(t)$. Moreover, each coefficient of $\Delta_K(t)$ is a finite type invariant, a linear combination and a product of finite type invariants is also a finite type invariant. Thus $a_n(K)$ is of finite type.
\end{proof}

We finish by a remark.
\begin{remark}
Let $K$ be the $(2, q)$ torus knot. Each component of the non--abelian part of the character variety  is parametrized by $\rho \mapsto \mathrm{tr} \rho(\mu)$. The invariant $Z_N(K)$ can be written as a single integral on the $\SU$-moduli space of the knot group as follows. 

Let $\widehat{R}(\mathscr{M}(2,q)) = \mathrm{Hom}^{\mathrm{irr}}(\Pi(2,q); \SU)/\SO$ be the $\SU$-moduli space and let $${CS}_{\mathscr{T}(2,q)} \colon \widehat{R}(\mathscr{M}(2,q)) \to S^1$$ denote the map defined by ${CS}_{\mathscr{T}(2,q)} (\rho) = CS_\ell$ if $\rho$ lies in the $\ell$-th component of $\widehat{R}(\mathscr{M}(2,q))$. Recall that $\widehat{R}(\mathscr{M}(2,q))$ is a $1$-dimensional smooth manifold (see~\cite{JDFourier}). One can prove the following formula:
\begin{equation}\label{ZWittSU}
Z_N({\mathscr{T}(2,q)}) = \frac{1}{\pi} \int_{\widehat{R}(\mathscr{M}(2,q))} \frac{{{CS}_{\mathscr{T}(2,q)}}^{-N}}{\sqrt{|\mathbb{T}_\lambda^{\mathscr{T}(2, q)}|}} \, \omega^{\mathscr{T}(2,q)}.
\end{equation}
Here $\omega^{\mathscr{T}(2,q)}$ denotes the $1$-volume form on $\widehat{R}(\mathscr{M}(2,q))$ defined in~\cite[\S~6]{JDFourier} using the Reidemeister torsion. Formula~(\ref{ZWittSU}) can be considered as an analogue for knots of the Witten formal integral for closed $3$-dimensional manifolds.
\end{remark}



%


\begin{thebibliography}{99}
%
\bibitem{Burde}
G. Burde, \emph{Darstellungen von Knotengruppen} (German), Math. Ann. \textbf{173} (1967) 24--33.
%
\bibitem{CS:1983}
M.~Culler and P.~Shalen, \emph{Varieties of group representations and
  splittings of $3$-manifolds}, Ann. of Math. \textbf{117} (1983) 109--146.
%
\bibitem{DeRham}
G. de Rham, \emph{Introduction aux polyn™mes d'un n\oe ud} (French), Enseignement Math. (2) \textbf{13} (1967) 187--194.
%
\bibitem{JDFourier}
J. Dubois, \emph{Non abelian Reidemeister torsion and volume form on the $\SU$-representation space of knot groups}, {Ann. Institut Fourier}, \textbf{55} (2005) 1685--1734.
%
\bibitem{JDFibre}
J. Dubois, \emph{Non abelian twisted Reidemeister torsion for fibered knots}, to appear in \emph{Canad. Math. Bull.} {arXiv:math.GT/0403304}.
%
\bibitem{Hikami}
K. Hikami, \emph{Volume Conjecture and asymptotic expansion of $q$-series}, Experiment. Math., \textbf{12} (2003) 319--337.
%
\bibitem{HK}
K. Hikami and A.N. Kirillov, \emph{Torus knot and minimal model}, Phys. Lett. B \textbf{575} (2003) 343--348.
%
\bibitem{Kas95}
R. Kashaev, \emph{A link invariant from dilogarithm}, Mod. Phys. Lett. A \textbf{10} (1995) 1409--1418.
\bibitem{Kas97}
R. Kashaev, \emph{The hyperbolic volume of knots from quantum dilogarithm}, Lett. Math. Phys. \textbf{39} (1997) 269--275.
\bibitem{KasTir}
R. Kashaev and O. Tirkkonen, \emph{Proof of the volume conjecture for torus knots}, J. Math. Sci. New York \textbf{115} (2003) 2033--2036.
%
\bibitem{KKK} 
P. Kirk and E. Klassen, \emph{Chern--Simons invariants of 3-manifolds decomposed along tori and the circle bundle over the representation space of $T^2$}, Commun. Math. Phys. \textbf{153} (1993) 521--557.
%
%
\bibitem{Klassen:1991}
E. Klassen, \emph{Representations of knot groups in $\SU$,} Trans. Amer. Math. Soc. \textbf{326} (1991) 795--828.
%
\bibitem{HPS}
M. Heusener, J. Porti and E. Su\'arez, \emph{Deformations of reducible representations of $3$-manifold groups into $\SL$}, J. reine angew. Math. \textbf{530} (2001) 191--227.
%
\bibitem{Le}
T. Le,  \emph{Varieties of representations and their subvarieties of cohomology jumps for certain knot groups}, Russian Acad. Sci. Sb. Math.\textbf{78} (1994) 187--209.
%
\bibitem{Milnor:1962}
J. Milnor, \emph{A Duality Theorem for Reidemeister Torsion}, Ann. of Math. \textbf{76} (1962) 134--147.
%
\bibitem{Milnor:1966}
J. Milnor, \emph{Whitehead torsion}, Bull. Amer. Math. Soc. \textbf{72} (1966) 358--426.
%
\bibitem{Mur}
H. Murakami and J. Murakami, \emph{The colored Jones polynomials and the simplicial volume of a knot}, Acta Math. \textbf{186} (2001) 85--104.
%
\bibitem{Mur2}
H. Murakami, \emph{Asymptotic behaviors of the colored Jones polynomials of a torus knot}, Internat. J. Math. \textbf{15} (2004) 547--555.
%
\bibitem{Ohtsuki}
T. Ohtsuki, {Problems on Invariants of Knots and $3$-Manifolds}, \emph{Geom. Topol. Monogr.} \textbf{4} (2002) 377--572.
%
\bibitem{Park:1997}
J. Park, \emph{Half-density volumes of representation spaces of some $3$-manifolds and their application}, Duke Math. J. \textbf{86} (1997) 493--515.
%
\bibitem{Porti:1997}
J.~Porti, \emph{{Torsion de Reidemeister pour les vari\'et\'es hyperboliques}},
  vol. 612, Mem. Amer. Math. Soc., 1997.
%
\bibitem{Riley}
R. Riley, \emph{{Nonabelian representations of $2$-bridge knot groups}}, Quart. J. Math. Oxf. \textbf{35} (1984) 191--208.
%
\bibitem{RSW}
T. Ramadas, I. Singer, and J. Weitsman, \emph{Some comments on Chern--Simons gauge theory},  Comm. Math. Phys. \textbf{126} (1989) 409--420. 
%
\bibitem{Turaev:1986}
V. Turaev, \emph{Reidemeister torsion in knot theory}, \emph{English version, } Russian Math. Surveys \textbf{41} (1986) 119--182.
%
\bibitem{Turaev:2000}
V. Turaev, \emph{Introduction to combinatorial torsions}, Birkh\"auser 2001.
%
\bibitem{Turaev:2002}
V. Turaev, \emph{Torsions of $3$-dimensional manifolds}, {Birkh\"auser} 2002.
%
\bibitem{YY2}
Y. Yamaguchi, \emph{The limit values of the non--abelian twisted Reidemeister torsion associated to knots,} preprint arXiv:math.GT/0512277 (2005)  
%
\bibitem{Zheng} 
H. Zheng, \emph{Proof of the volume conjecture for Whitehead double of tours knots}, preprint available at arXiv:math.GT/0508138.
\end{thebibliography}
\end{document}